\documentclass[12pt]{amsart}
\usepackage{xcolor}
\usepackage{graphicx}
\usepackage{graphics}
\usepackage{amsmath}
\usepackage{amscd}
\usepackage{latexsym}
\usepackage{subfigure}
\usepackage{caption}
\usepackage{hyperref}
\usepackage{dsfont}
\usepackage{amsthm}
\usepackage{amssymb}
\usepackage{placeins}{\Large}

\begin{document}

\textwidth 5.9in
\textheight 7.9in

\evensidemargin .75in
\oddsidemargin .75in

\newtheorem{Thm}{Theorem}
\newtheorem{Lem}[Thm]{Lemma}
\newtheorem{Cor}[Thm]{Corollary}
\newtheorem{Prop}[Thm]{Proposition}
\newtheorem{Rm}{Remark}

\def\a{{\mathbb a}}
\def\C{{\mathbb C}}
\def\N{{\mathbb N}}
\def\A{{\mathbb A}}
\def\B{{\mathbb B}}
\def\D{{\mathbb D}}
\def\E{{\mathbb E}}
\def\R{{\mathbb R}}
\def\P{{\mathbb P}}
\def\S{{\mathbb S}}
\def\Z{{\mathbb Z}}
\def\O{{\mathbb O}}
\def\H{{\mathbb H}}
\def\V{{\mathbb V}}
\def\Q{{\mathbb Q}}
\def\Cn{${\mathcal C}_n$}
\def\CM{\mathcal M}
\def\CG{\mathcal G}
\def\CH{\mathcal H}
\def\CT{\mathcal T}
\def\CF{\mathcal F}
\def\CA{\mathcal A}
\def\CB{\mathcal B}
\def\CD{\mathcal D}
\def\CP{\mathcal P}
\def\CS{\mathcal S}
\def\CZ{\mathcal Z}
\def\CE{\mathcal E}
\def\CL{\mathcal L}
\def\CV{\mathcal V}
\def\CW{\mathcal W}
\def\IC{\mathbb C}
\def\IF{\mathbb F}
\def\IK{\mathcal K}
\def\IL{\mathcal L}
\def\IP{\bf P}
\def\IR{\mathbb R}
\def\IZ{\mathbb Z}

\title{On the Smale Conjecture for Diff$(S^4)$}
\author{Selman Akbulut}

\address{Gokova Geometry Topology Institute,
Akyaka, Mehmet Gokovali Sokak,
No:53, Ula, Mugla, Turkey}
\email{akbulut.selman@gmail.com}

\date{\today}
\begin{abstract} 
Recently Watanabe disproved the Smale Conjecture for $S^4$, by showing Diff$(S^{4})\neq SO(5)$. He showed this by proving that their higher homotopy groups are different. Here we prove this more directly by showing $\pi_{0}$Diff$(S^{4})\neq 0$, otherwise a certain loose-cork could not possibly be a loose-cork. 
\end{abstract}
 \date{}
\maketitle

\setcounter{section}{-1}

\vspace{-.3in}

\section{An exotic diffeomorphism }\label{exoticdiff}

Here we prove  $\pi_{0}$Diff$(S^{4})\neq 0$, by showing that if this is not true, then the loose-cork defined in \cite{a2} could not possibly be a loose-cork. The group $\pi_{0}$Diff$(S^{4})= \pi_{0}$Diff$(B^{4}, S^{3})$ can be calculated from the homotopy exact sequence of the following fibration (where $I=[0, 1]$, $\dot{I}=\partial I$)
\begin{equation}\mbox{Diff}(S^{3}\times I, S^{3}\times \dot{I}) \to \mbox{Diff}(B^{4},S^{3}) \to
\mbox{Emb}(B^{4}, \mbox{Int}B^4)\end{equation} as the free part part of  $\pi_{0}\mbox{Diff}(S^{3}\times I, S^{3}\times \dot{I}) =\pi_{0}\mbox{Diff}(B^4, S^{3}) \oplus Z_{2}$. $\mbox{Emb}(B^{4}, \mbox{Int}B^4)$ is path connected, and the $\Z_{2}$ summand comes from Dehn twisting $S^{3}\times I=B^{4}-\mbox{Int}(B^{4})$ along $S^{3}$ by using $\pi_{1}SO(4)$.  

\vspace{.1in}

A  nontrivial element of  $\pi_{0}$Diff$(B^{4}, S^{3})$ can be described as follows:  Let $T \subset S^{3}$ be a tubular neighborhood of the figure-8 knot $K$ lying inside of the $3$-ball $B^{3}$. Call $C=T\times J\approx S^{1}\times B^{2}\times J\subset B^{3} \times J \subset S^{3}\times J$, where $J =[-1, 1]$. Define a diffeomorphism $f_{n}: C \to C$ (where $n \in \N$)  
$$f_{n}(x,y, t )= \left(xe^{-2\pi i n^{2}t^{2}},y, t \right)$$


\noindent As $t$ runs from $-1$ to $1$, $f_{n}$ rotates  $T$ forward (and $f^{-1}_{n}$ backwards) Let $\bar{C}=C\cup \; S^{3}\times N(\dot{J})\subset S^{3}\times J $, where $N(\dot{J})$ is a small closed neighborhood of the boundary $\dot{J}\subset J$. After resizing $J$ and extending domain of definition we extend $f_{n}: \bar{C}\to \bar{C}$.

   \begin{figure}[ht]  \begin{center} \includegraphics[width=.16\textwidth]{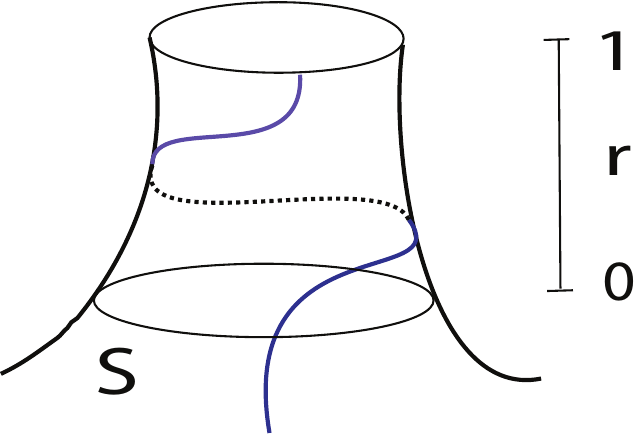}       
\caption{Dehn twist}      \label{a2} 
\end{center}
 \end{figure}
 
$K$ is a fibered knot with fibration $S^{3} - T\to S^1$, and the fiber punctured torus $S$. Let $\delta_{s}: S\to S$ be Dehn twisting diffeomorphism along the boundary parallel curve. That is, if we identify a collar $N$  of $\partial S$ with $N\approx S^{1}\times (0,1] $, then $\delta_{s}: N \to N$, is 
given by $\delta_{s}(x,r)=(xe^{2\pi is r} ,r)$. Since $\delta_{s}$ commutes with monodromy of the fibration $S^{3}-T \to S^{1}$, it induces a diffeomorphism $\bar{\delta}_{s}: S^{3}-T \to S^{3}-T$, which is $e^{2\pi is}$ rotation along  $\partial T$ in $K$ direction. Clearly $\bar{\delta}_{s}$ is supported near $T$. 

\vspace{.1in}

 Now define  $g_{n}: (S^{3} -T) \times J  \to  (S^{3} -T)\times J$ by  $g_{n}(x,t)= (\bar{\delta}_{s}(x), t)$ where $s=t^{2}n^{2}$. By considering orientations, we can decompose 
 $$ S^{4}= (S^{3} -T)\times J \smile - (T\times J ) $$ and  extend $g_{n}$ to a diffeomorphim $\phi_{n}: S^{3}\times J \to S^{3}\times J$ by letting it to be $f^{-1}_{n}$ on $-(T \times J)$ (accounting the change of orientation in gluing).
 So, as t runs from $-1$ to $1$, at each $t$ level $ \phi_n$ simultaneously Dehn twists the pages  of the fibration $S^{3}-T \to S^{1}$ by $2\pi n^{2} t^{2}$ amount to left, while rotating $T$ in the longtitutional direction $-2\pi n^{2}t^{2}$ to the right.

\vspace{.05in}

Since $\phi_{n}$ is supported near $C$, it fixes some small neighborhood of a vertical arc $p \times J$, where $p\in S^{3}-T$.
  Let $N(p) \approx B^3$ be a small  ball neighborhood of $p$, then by the identifications $B^{4} \approx B^{3}\times J  \approx ( S^{3} - N(p) )\times J $. Then by restriction, we get a diffeomorphism 
  $$\phi_{n}^{0}: (B^{4}, S^{3}) \to (B^{4}, S^{3})$$

 \begin{figure}[ht]  \begin{center}
 \includegraphics[width=.52\textwidth]{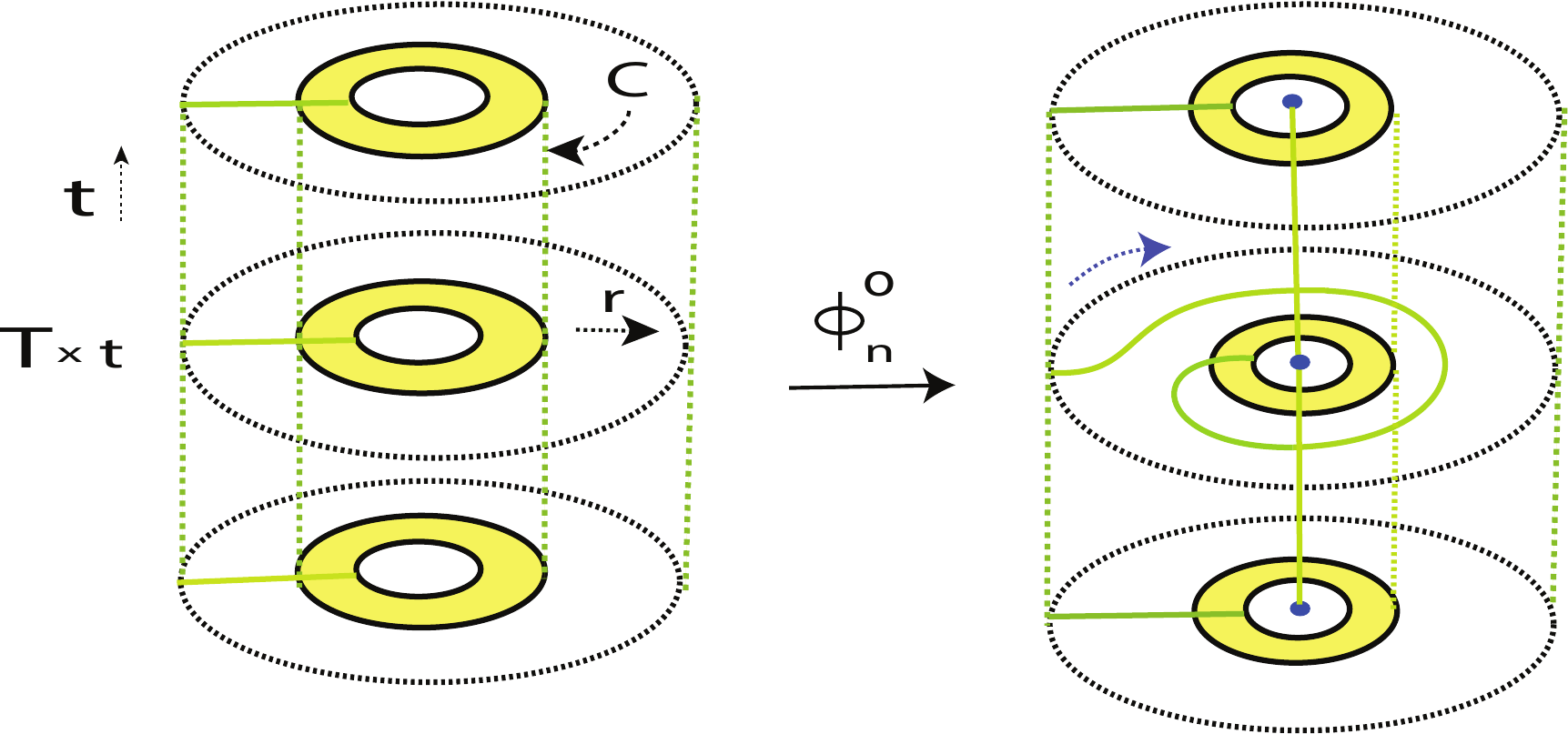}       
\caption{$\phi_{n}^{0} : B^{3}\times J \to B^{3}\times J$}     \label{s81} 
\end{center}
 \end{figure}

    The vertical outer boundary of the wall of Figure~\ref{s81} corresponds to the vertical outer boundary of $N(p) \times J$, where $\phi_{n}^{0}$ is identity. We will call $\phi^{0}_{n}$ {\it Dehn twisting diffeomorphism of $B^{4}$ along $C$}. More generally, when $K\subset S^{3}$ is any fibered knot and $T^{2}$ is the boundary of its tubular neighborhood, we will also refer the resulting $\phi_{n}^{0}: B^{3}\times J \to  B^{3}\times J $ Dehn twisting   of $B^{4}$ along  $T^{2}$.  We will show that $\phi_{n}^{0}$ is not isotopic to identity fixing the boundary, but if we allow it to move the boundary it is isotopic to identity by the so called ``swallow-follow isotopy".
    
    \newpage
          
  \begin{Thm}
 The diffeomorphisms $\phi_{n}^{0}:(B^4, S^3)\to (B^4, S^3)$ are not isotopic to each other fixing the boundary, for distinct integers $n>0$. \end{Thm}

 \proof  
 We will prove this theorem by constructing a properly imbedded disk $D\subset B^{4}$ which can not be isotopic to $\phi_{n}^{0}(D)$ rel boundary, otherwise certain associated cork (which we will define in the next section) would not be a cork: Let $K_{0}$ be the punctured $K$, and let $D = K_{0}\times J$  be the obvious disk in $B^{4}=B^{3}\times J$, which $K\#-K$ bounds.

 \vspace{.05in}
 
\begin{figure}[ht]  \begin{center}
\includegraphics[width=.75\textwidth]{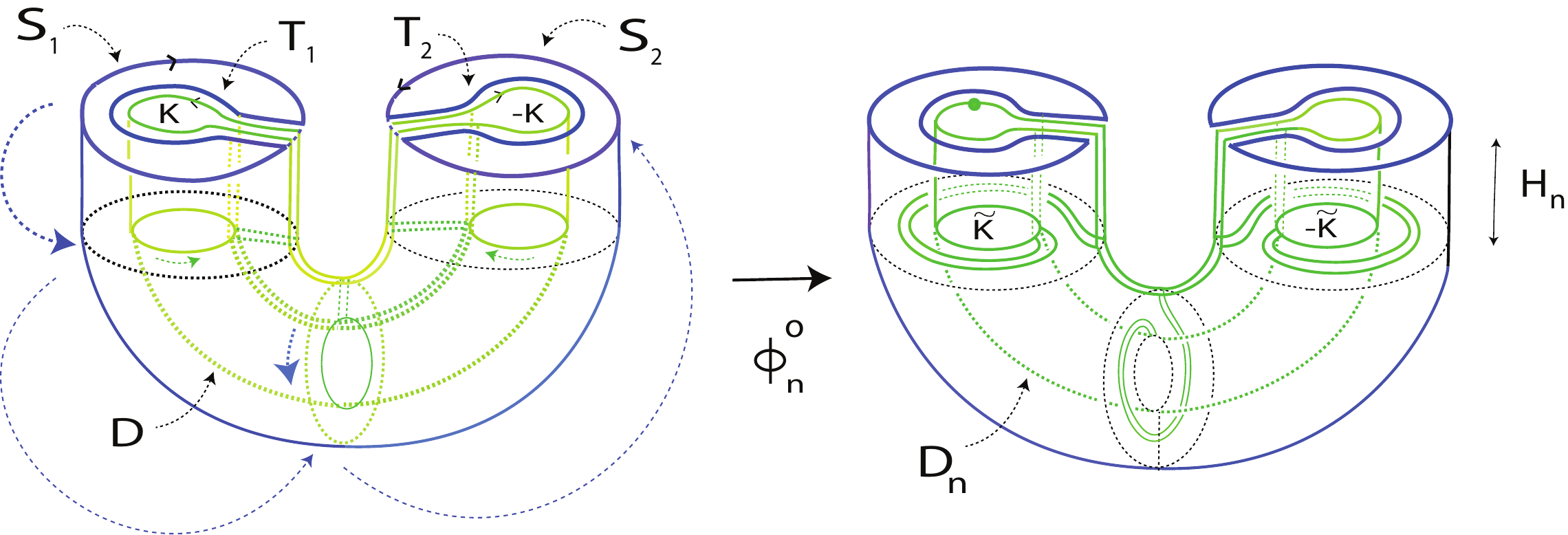}       
\caption{$\phi^{o}_{n} : B^{3}\times J \to B^{3}\times J$} \label{s8c} 
\end{center}
 \end{figure}

 Figure~\ref{s8c} is the same as Figure~\ref{s81}, drawn by bending the figure. Then take $D_{n}$ be the image of this disk $D$ under $\phi_{n}$. $D_{n}$ is the union of a concordance $H_{n}$ from $K\# -K$ to $\bar{K}\# -\bar{K}$ induced by $\phi_{n}$, and the disk $\bar{D}_{n}=\bar{K}_{0} \times I $, where $\bar{K}$ is the other end of the concordance $H_{n}$ as shown in the second picture of Figure~\ref{s8c}. So we have $\phi^{o}_{n}:(B^{4},D) \to (B^{4},D_{n}) $ 

\begin{figure}[ht]  \begin{center}
 \includegraphics[width=.6\textwidth]{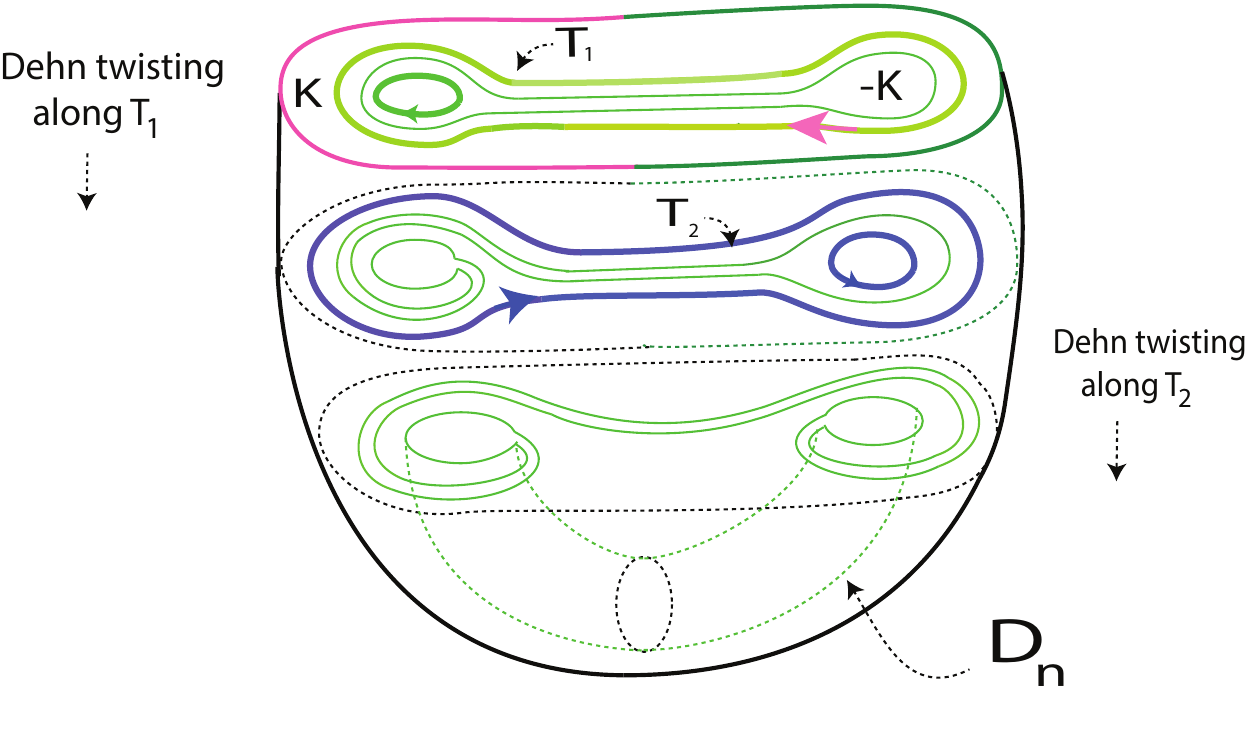}       
\caption{Forming $D_{n}$ by swallow-follow isotopies}  \label{s8k} 
\end{center}
 \end{figure}

\begin{figure}[ht]  \begin{center}
 \includegraphics[width=.7\textwidth]{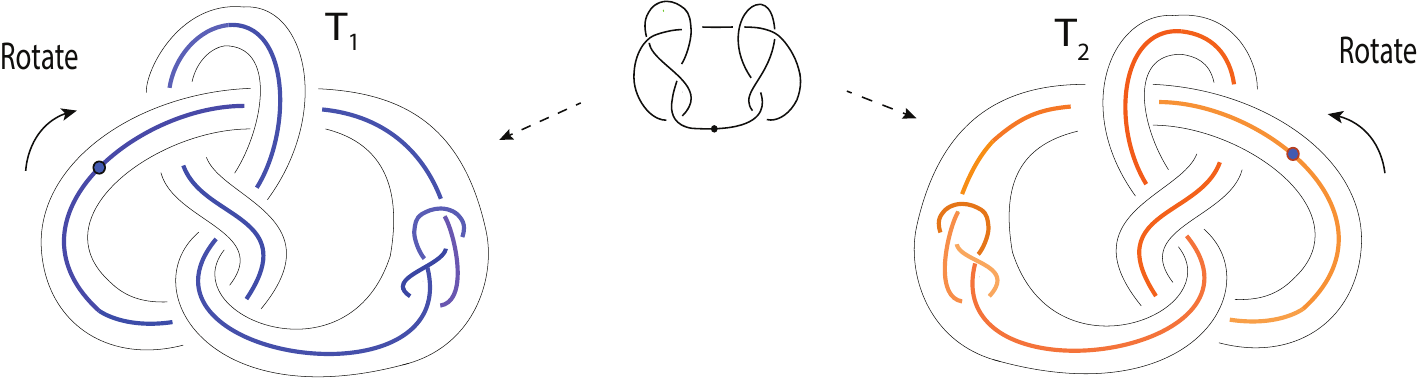}       
\caption{Swallow-follow isotopies of $K\#-K$} \label{s9b} 
\end{center}
 \end{figure}

 We shell see, although $D$ and $D_{n}$ are not isotopic fixing their boundaries on $\partial B^{4}$, they are isotopic moving the boundaries, through an isotopy which is  called  {\it `swallow-follow'}. This isotopy is a composition of two isotopies of $K\#-K$: By first moving $-K$ along $K$ in its tubular neighborhood, then moving $K$ along $-K$ in its tubular neighborhood, than capping the resulting $K\#-K$ by the obvious disk it bounds in $B^{4}$. Various pictures of these moves are described in Figures~\ref{s8k}, ~\ref{s9b}). 
 
 \vspace{.1in} 
 
 We can also explain swallow-follow isotopy from Figure~ \ref{s8c}, by first constructing two disjoint $2$-toruses $T_1$ and $T_2$ in  $\partial (B^{3}\times J)$, enclosing $K$ and $-K$ respectively, not intersecting $D$. This can be done by piping  boundaries of the tubular neighborhoods  $T_{1}'$, $T_{2}'$, of $K$ and $-K$ to large disjoint spheres $S_1$ and $S_2$, enclosing them. That is $T_{i}=T_{i}'\#S_{i}$, with $i=1,2$ as shown in  Figure~ \ref{s8c}. There is an isotopy $S_{1} \rightsquigarrow S_{2}$, indicated by the arrows of Figure~ \ref{s8c}. This has an affect of replacing the core $K$ in its tubular neighborhood by $K\#-K$, and then replacing the core $-K$ in its tubular neighborhood by $-K\#K$.

\vspace{.1in}

 $T_{1}$ and $T_{2}$ are the two ends of an imbedded copy of $T\times J \subset B^{3}\times J$, which is disjoint from $D$. Here we can take this $T$ to be solid torus. Dehn twisting diffeomorphism of $B^{3}\times J$ along $T\times J$ takes $D$ to $D_{n}$. This induces diffeomorphism  between  $B^{3}\times J$ with $D$ removed (carved) and $B^{3}\times J$ with $D_{n}$ removed. As can be seen from Figure~\ref{s8k}, this diffeomorphism amounths to gluing two isotopies, each taking place in tubular neighborhoods of  $K$ and $-K$ respectively. The effect of this is to tie a small copy of $-K$ to $K$, and move along $K$, then tie a small copy of $K$ to $-K$ and move along $-K$.
 
   \vspace{.05in}

 \newpage

  \section{Forming the infinite order cork }\label{infcork}

\begin{figure}[ht]  \begin{center}
 \includegraphics[width=.23\textwidth]{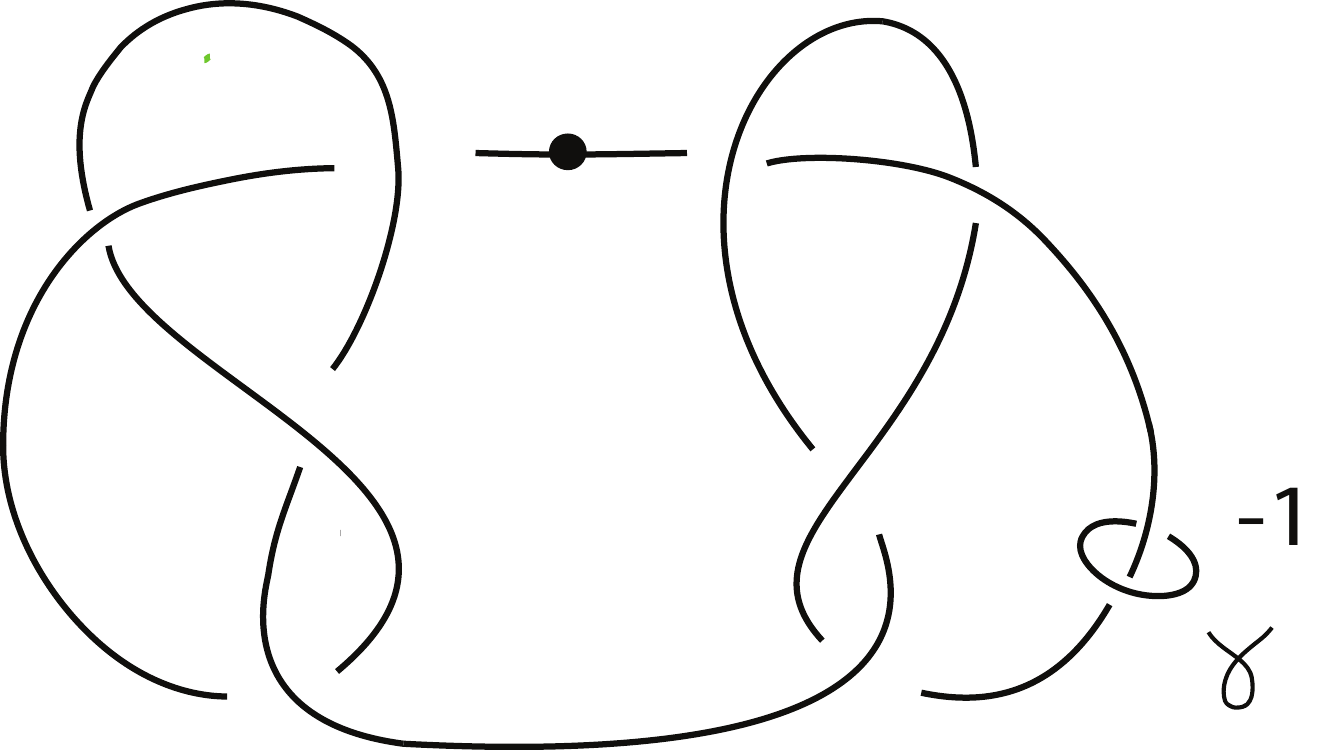}       
\caption{W}      \label{s3} 
\end{center}
 \end{figure}

 Next, from $\phi_{n}$ we form an infinite order cork $(W, \tau_{n})$: We do this by removing the slice disk $D$ from $B^4$ then attaching $-1$ framed $2$-handle to the meridian linking circle $\gamma$ (Figure~\ref{s3}), i.e. $W^{*}: = B^{4} - N(D)$ and
 
 $$W=W^{*} + h^2_{\gamma}$$ 
 
\begin{figure}[ht]  \begin{center}
 \includegraphics[width=.55\textwidth]{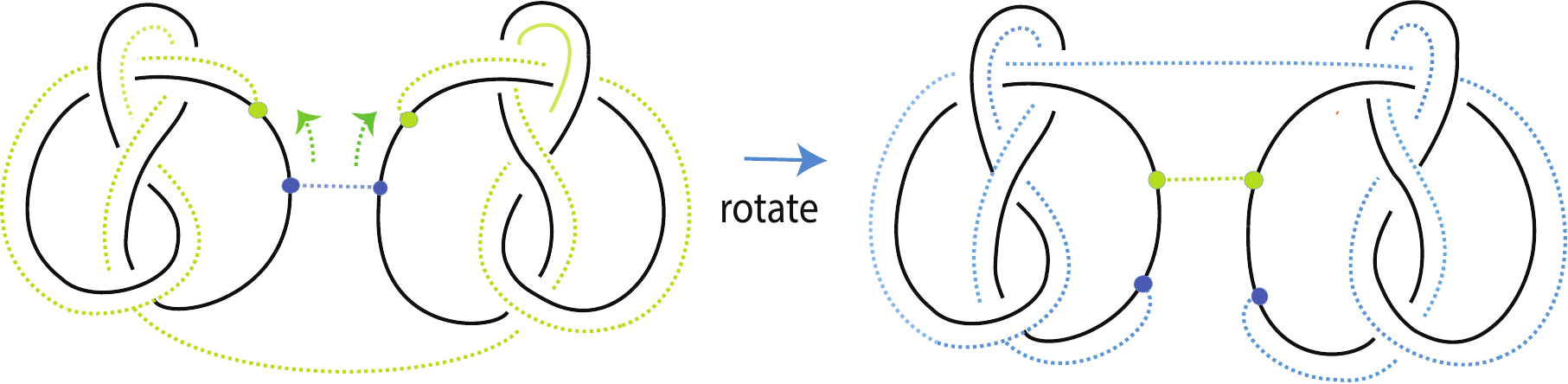}       
\caption{Rotating $K\#-K$}     \label{r3} 
\end{center}
 \end{figure}

\begin{figure}[ht]  \begin{center}
 \includegraphics[width=.55\textwidth]{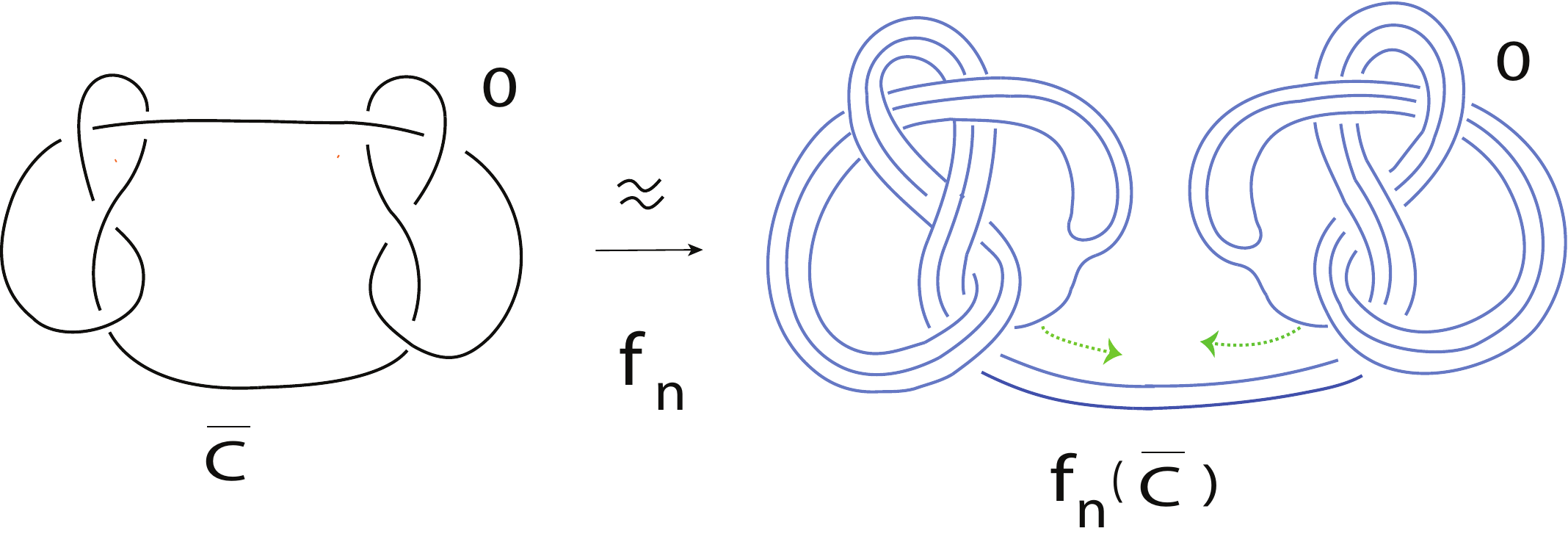}       
\caption{ }      \label{r1} 
\end{center}
 \end{figure}

\begin{figure}[ht]  \begin{center}
 \includegraphics[width=.55\textwidth]{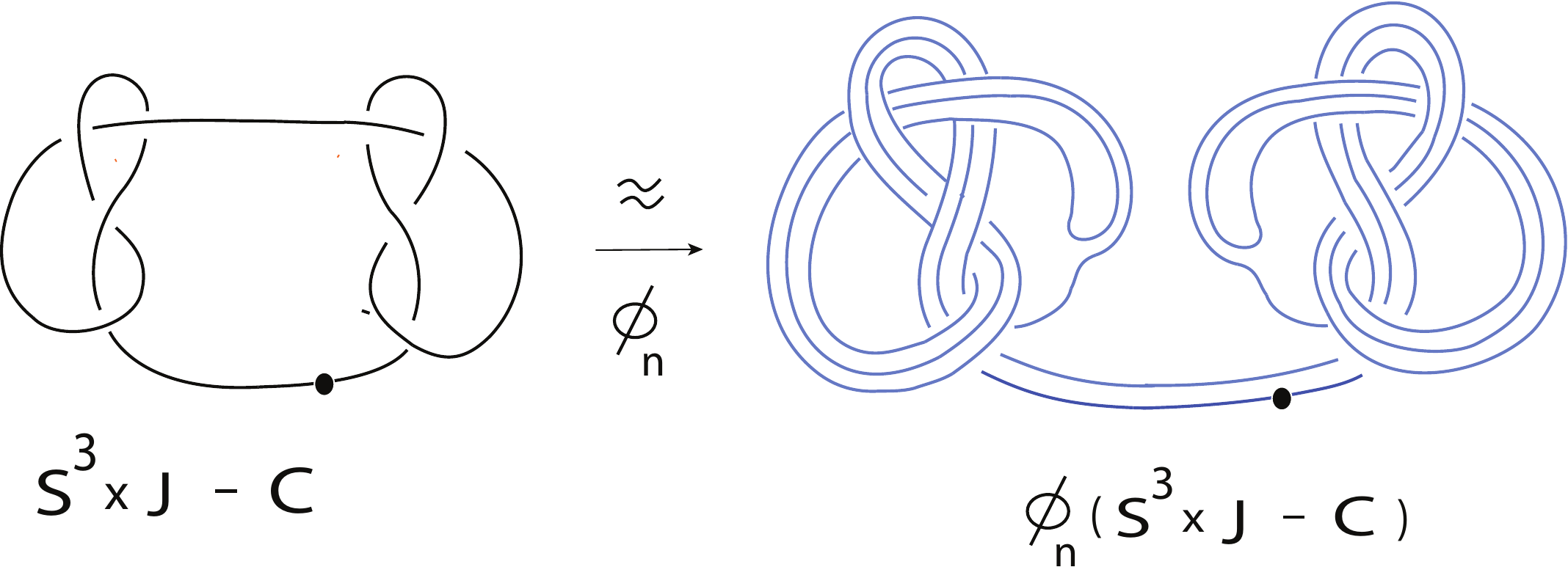}       
\caption{ }      \label{r2} 
\end{center}
 \end{figure}

The cork twisting map $\tau_{n}: \partial W \to \partial W$ is given by $f_{n}$ (Figures~ \ref{r1} \ref{r2}) which is induced from the diffeomorphism $\phi_{n}: W^{*} \to  W_{n}^{*}$, which keeps $\gamma$ fixed and twists $\partial W^{*}$ by $f_{n}$. When viewed $\phi_{n}$ as a self map of $W^{*}$, $f_{n}$ becomes visible as a combination of two Dehn twists in 
$\partial W^{*}$ opposite direction. This appears as a rotation along the connecting torus in the middle of Figure~\ref{s8c} as indicated in Figure~\ref{r3} (recall that $\partial W^{*}$ consists of two copies of $S^{3} - N(K)$ glued along their boundary toruses \cite{a1}). This rotation $\tau_{n}$ on the boundary can not extend inside $W$. 

 Next we will describe the images of the handles of $S^{3}\times J$ by $\phi_{n}$. Figure~\ref{r1} describes $\bar{C}$ and $f_{n}(\bar{C})$. Here $\bar{C}$  is drawn as a round $2$-handle attachment to two disjoint copies of $S^{3}\times I$ ($0$- and $4$- handles are not drawn since hey are attached uniquely). Handlebody of $f_{n}(\bar{C})$ is constructed similarly as a round $2$-handle attachment, except in this case while attaching the round $2$-handle, we rotate $T\times J$ in $J$ direction.  
 
 \noindent That is, rotate $n$ and $-n$-times near each of its boundary components. Recall that, $W_{n}$ is carved from $S^{3}\times J$ by using $f_{n}(C)$.
  
\newpage

 Figure~\ref{r2} describes the complement $S^{3}\times J - C \approx (S^{3}-T) \times J$,  and its image $g_{n} (S^{3}\times J - C)$. The left picture of Figure~\ref{r2} is  the complement of the left picture of Figure~\ref{r1} in $S^{3}\times J$. The right picture of  Figure~\ref{r2} is the complement of the right picture of Figure~\ref{r1} (recall decomposing $S^{4}$ using a slice knot, by attaching and carving $2$-handles from the two hemispheres $B^{4}_{\pm}$  as in 14.3 of \cite{a1}). Here $\phi_{n}$ is the important part, which is used in the construction of the cork.
 
 \vspace{.1in}
 
 The two right pieces are the nontrivial diffeomorphism induced by the ``swallow-follow" isotopy. $\phi_{n}|_{C}$ is longitudinal rotation of $T$ along $J$, and $\;g_{n}(S^{3}\times J -C)$ describes a carved out disk  $D_{n}\subset B^{4}$ bounded by $K\#-K$, which is the disk obtained by concatenating a swallow-follow concordance of $K\#-K$, followed by the standard disk $D$ which $K\#-K$ bounds in $B^{4}$ and  
$D_{n}=\phi_{n}(D)$.

\section{Constructions and proofs }\label{introduction}

\begin{figure}[ht]  \begin{center}
 \includegraphics[width=.55\textwidth]{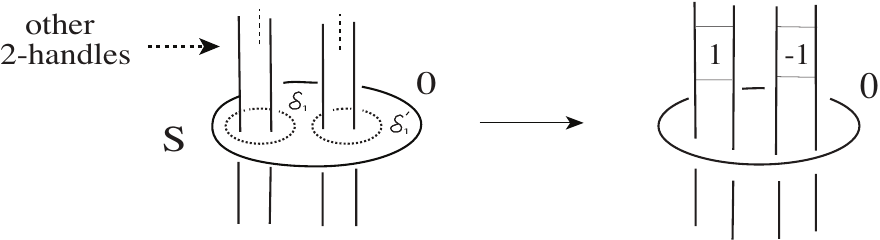}       
\caption{Gluck twisting to $M$ along $S$}      \label{s1} 
\end{center}
 \end{figure}

Let $M^4$ be a smooth $4$-manifold, and $S^2 \times D^2\subset M^4$ be the tubular neighborhood of an imbedded $2$-sphere $S$, denoted as $0$-framed circle in the left picture of  Figure \ref{s1}.  Gluck twisting of $M$ along $S$ is the operation of cutting out $S^2\times D^2$ from $M$, and regluing by the nontrivial diffeomorphism of $S^2 \times S^1$. The affect of this operation on the handlebody is indicated in Figure~\ref{s1} (\cite{a1}). This is the operation of arbitrarily separating $2$-handle strands going through $S$ into two groups, and applying $1$ twist across one group and $-1$ twist across the other.  If $H\subset M$ is an imbedded cylinder with boundary components $\delta_1$, $\delta_1'$ away from the $2$-handles, this operation corresponds to twisting $M$ along $H$. 

\begin{figure}[ht]  \begin{center}
\includegraphics[width=.6\textwidth]{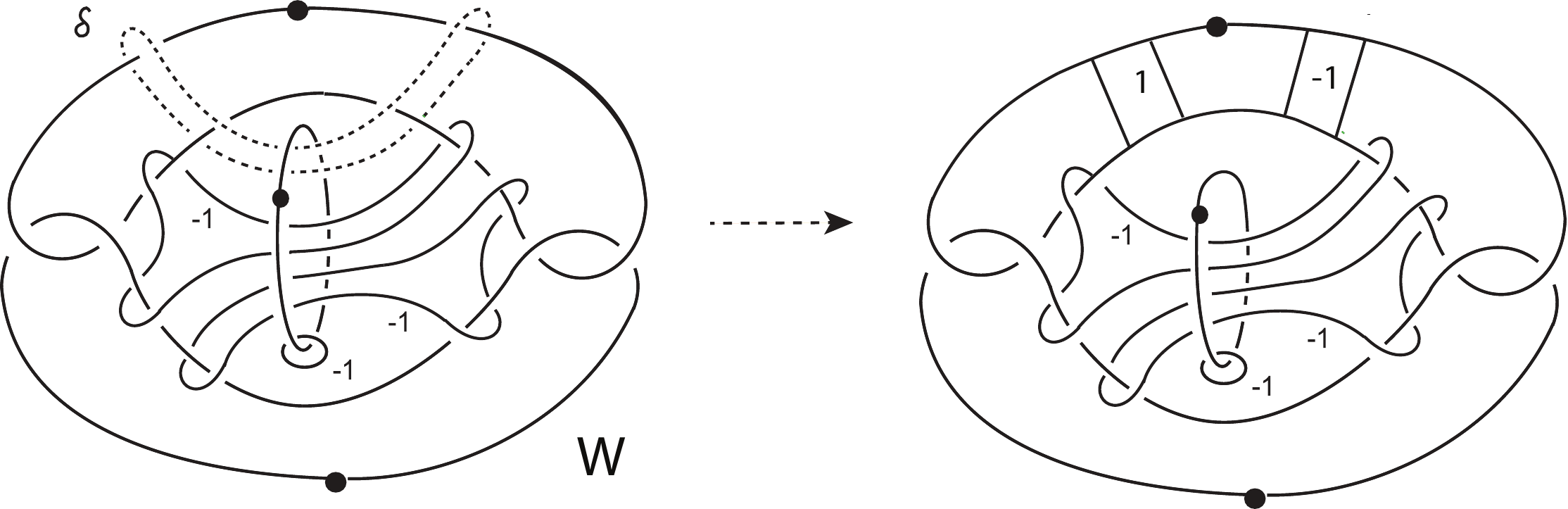}       
\caption{$\delta$-move diffeomorphism  $W\approx W$}   \label{s2} 
\end{center}
 \end{figure}

Now recall the infinite order loose-cork $(W,h)$, defined in \cite{a2} and \cite{g}. As discussed in \cite{a2}, $W$ is the contractible manifold of Figure~\ref{s2}, and the cork automorphism $h:  \partial W \to \partial W$ is given by the ``$\delta$-move'', which is indicated  in Figure~\ref{s2}. This operation is similar to Gluck twisting, where $S$ is replaced by the unknotted circle $\delta$. Main difference is, here we allow $1$-handles (circle-with-dots) go through the circle $\delta$.

 \begin{Rm}
 At first glance reader might think the right and left twists of $\delta$-move in Figure~\ref{s2}, would cancel each other and nothing happens. But this is not so, this induces a nontrivial diffeomoprphism of the boundary of $W$ (which is the cork automorphism). For example, the delta move in Figure 26 of \cite{a3}, alters the position of $\gamma$ in a nontrivial way. $\gamma$ could be the attaching circle of a $2$-handle on top of the cork.
\end{Rm}

 Figure~\ref{s3}  is an equivalent definition of $W$, which is the contractible manifold obtained by blowing down $B^{4}$ (6.2 of \cite{a1}) along the obvious ribbon disk $D\subset B^4$ bounding $K\#-K$, where $K$ is the figure $8$ knot. 

\vspace{.05in}

Figure~\ref{s3} version of $W$ can be obtained from Figure~\ref{s2} by ignoring the middle dotted circle ($1$-handle) of the left picture of Figure~\ref{s2}, then canceling other two  circles-with-dots with the $-1$ framed $2$-handles. This process takes the middle 
$1$-handle circle to the dotted $K\#-K$.

\vspace{.1in}

\begin{figure}[ht]  \begin{center}
\includegraphics[width=.65\textwidth]{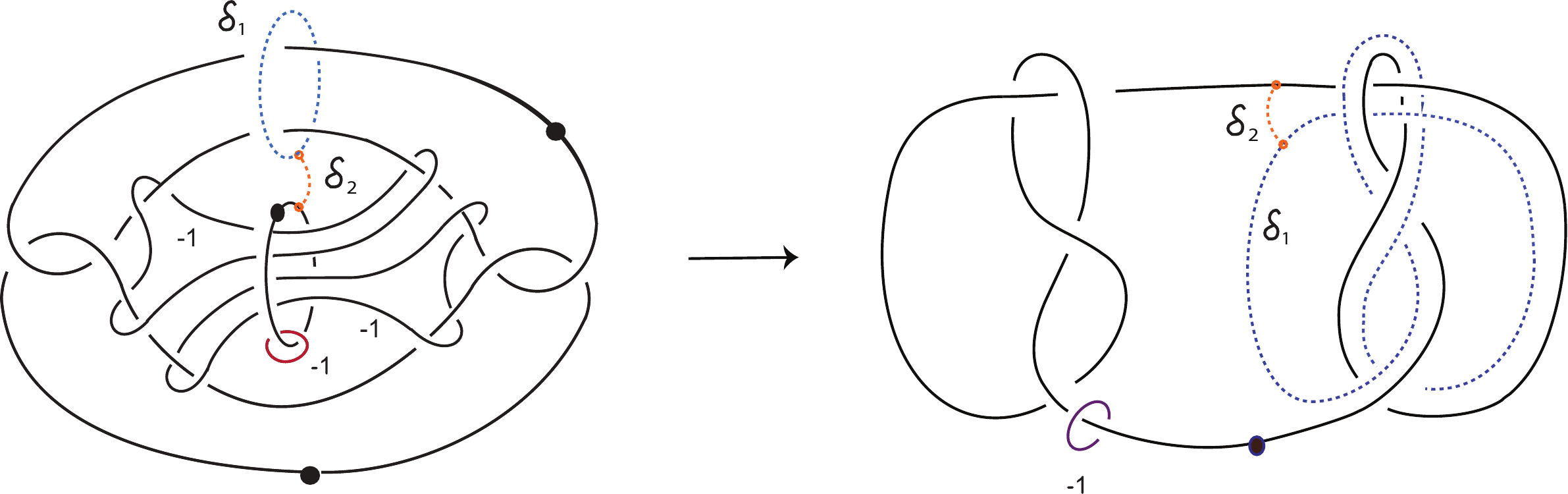}       
\caption{Two descriptions of W}      \label{s4} 
\end{center}
 \end{figure}

\begin{figure}[ht]  \begin{center}
\includegraphics[width=.65\textwidth]{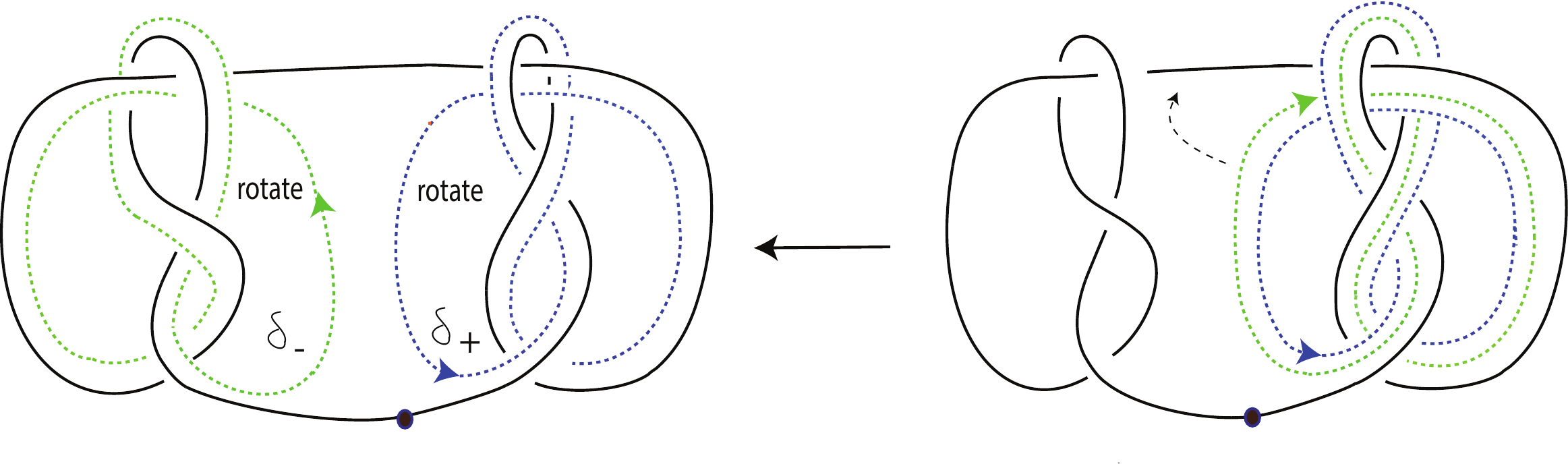}       
\caption{$\delta$-move}      \label{t1} 
\end{center}
 \end{figure}

   Figure~\ref{s4} indicates to where this process takes the curve $\delta_1$ from left to the right picture of this figure (this is a very crucial observation).  
   In the left picture  $\delta$-move corresponds  twisting along two parallel copies of  $\delta_{1}$ in the opposite direction. To see this in the right picture of Figure~\ref{s4}, we have to perform this operation along two parallel copies of the curve corresponding $\delta_{1}$, which is a  copy of figure-8 knot. One of the copies of the knot can be slid over the ribbon $1$-handle (shown in the figure) to put in the position of the left picture of Figure~\ref{t1}.

  So oppositely twisting W along $\delta_{-}$ and $\delta_{+}$ curves of Figure~\ref{t1} will result the desired $\delta$-move. This corresponds to altering the position of the standard ribbon disk $D\subset B^{4}$ which $K\#-K$ bounds to another ribbon disk $D_{n}$, which is obtained by concatenating the concordance induced by the isotopy with $D$. From our construction $D_{n}=\phi_{n}(D)$.  If $\phi_{n}$  was isotopic to identity relative to boundary, the infinite order cork automorphism $h: \partial W \to \partial W$  (as $n\to \infty$) would extend as a smooth diffeomorphism inside W, which is a contradiction.  \qed

  \begin{figure}[ht]  \begin{center}
 \includegraphics[width=.5\textwidth]{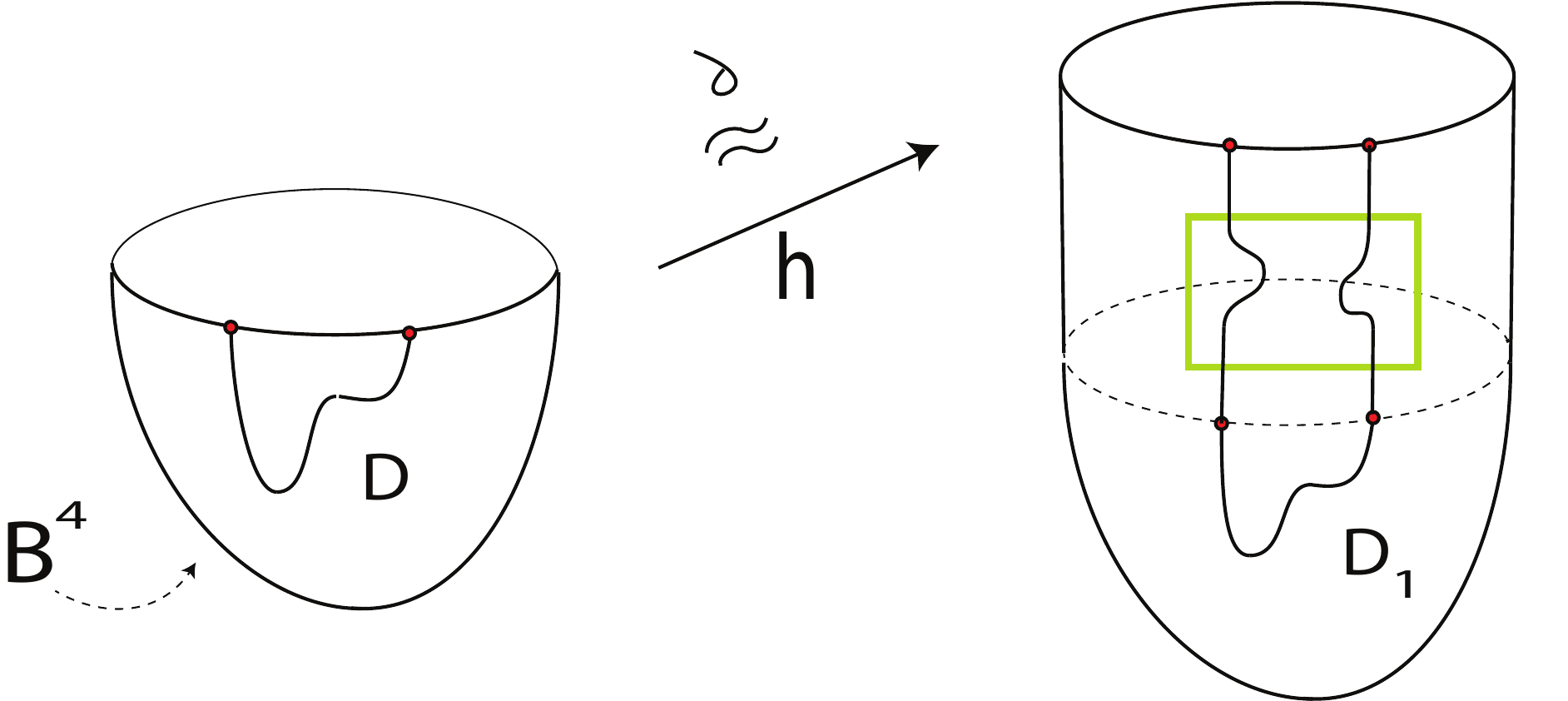}       
\caption{h}      \label{s6} 
\end{center}
 \end{figure}

 \begin{Rm}\label{boundary}
Diffeomorphism type of the carved-out $B^{4}$, carved along a properly imbedded disk $D\subset B^{4}$, does not change if we isotope $\partial D \subset S^{3}$ before carving, however blowing down isotopic disks $D, D'\subset B^{4}$, fixing $\partial D$ and $\partial D'$, can give corks  as in (\cite{a4}), and as in this example.
\end{Rm}

{\it Acknowledgements: I thank Michael Freedman for being a supportive friendly ear, and giving helpful suggestions while discussing this article. I would also like to thank Eylem Yildiz and Burak Ozbagci for helpful remarks; and thank Robion Kirby and Allen Hatcher for catching flaws in the previous versions of this article.}


\begin{thebibliography}{99999}

\bibitem[A1]{a1} S. Akbulut, {\em $4$-Manifolds}, Oxford Univ Press. ISBN-13: 978-0198784869

\bibitem[A2]{a2} S. Akbulut, {\em On infinite order corks}, PGGT, IP Paperback (2017) 151-157, \\  ISBN: 9781571463401,  \url{https://arxiv.org/pdf/1605.09348.pdf} 

\bibitem[A3]{a3} S. Akbulut, {\em Homotopy $4$-spheres associated to an infinite order loose cork}, \url{https://arxiv.org/pdf/1901.08299.pdf}

\bibitem[A4]{a4} S. Akbulut, {\em Corks and exotic ribbons in $B^4$}, European Journal of Math (2022)
\url{https://arxiv.org/pdf/2103.13967.pdf}
 
\bibitem[G]{g} R. E. Gompf, {\em Infinite order corks}, G$\&$T, vol.21, no.4 (2017)  2475-2484.

\bibitem[W]{w} T. Watanabe, {\em Some exotic nontrivial elements of the rational homotopy groups of Diff$(S^4)$}, \url{https://arxiv.org/pdf/1812.02448.pdf}

\end{thebibliography}
\end{document}